\documentclass[11pt,letterpaper,oneside]{amsart}

\usepackage{amsthm, amsmath, amssymb}
\usepackage{comment}
\usepackage{color}

\theoremstyle{plain}
\newtheorem{theorem}{Theorem}[section]
\newtheorem{lemma}[theorem]{Lemma}

\theoremstyle{definition}

\theoremstyle{definition}
\theoremstyle{remark}
\newtheorem{remark}{Remark}

 \newcommand{\be}{\begin{equation}}
 \newcommand{\ee}{\end{equation}}
 \newcommand{\bd}{\begin{displaymath}}
 \newcommand{\ed}{\end{displaymath}}
 \newcommand{\bea}{\begin{eqnarray}}
 \newcommand{\eea}{\end{eqnarray}}
 \newcommand{\beas}{\begin{eqnarray*}}
 \newcommand{\eeas}{\end{eqnarray*}}
 \newcommand{\bc}{\begin{center}}
 \newcommand{\ec}{\end{center}}

 \title{On  the uniqueness  of the  limiting solution to a strongly competing system }

\author{Avetik Arakelyan}
\address{Institute of Mathematics, NAS of Armenia, 0019 Yerevan, Armenia}
\email{arakelyanavetik@gmail.com}

\author{Farid Bozorgnia*}
\address{Department of Mathematics, Inst. Superior T\'{e}cnico, 1049-001 Lisbon, Portugal}
\email{bozorg@math.ist.utl.pt}

\thanks{A. Arakelyan was partially supported by State Committee of Science MES RA, in frame of the research project No.  16YR-1A017 }

\thanks{*The corresponding author,  F. Bozorgnia  was supported by the FCT post-doctoral fellowship
SFRH/BPD/33962/2009 }


\begin{document}

\begin{abstract}
	This work is devoted to prove uniqueness  result   for the positive solution to a  strongly competing system of Lotka-–Volterra type in the limiting configuration,  when the competition rate tends to infinity. Based on properties of limiting solution an alternative proof to show uniqueness is given.

		
\end{abstract}

\maketitle

\renewcommand{\theequation}{1.\arabic{equation}}
\setcounter{equation}{0}

\section{Introduction}


 Let $\Omega  \subset  \mathbb{R}^d$  be an open, bounded, and connected domain with smooth boundary. We take $m$  to be an integer number.  The aim of this paper is to investigate the  uniqueness  of solution for  a competition-–diffusion system of Lotka-–Volterra type, with Dirichlet boundary conditions as the competition rate tends to infinity.
This model of  strongly competing systems  have  been extensively studied  from different point  of views,    see   \cite{3,5, 6,7,8,9} and references therein.

The  model describes  the steady state   of $m$ competing species coexisting in the same area $\Omega.$
Let $u_{i}(x)$ denote the population density of the $i^{th}$ component.   The following system  shows the steady state of  interaction between $m$  components
\begin{align}\label{s1}
\begin{cases}
\Delta  u_{i}^{\varepsilon}=   \frac{ 1 }{\varepsilon}  u_{i}^{\varepsilon} \sum\limits_{j \neq i}   u_{j}^{\varepsilon} (x)\qquad\qquad & \text{ in  } \Omega,\\
u_{i}^{\varepsilon} \ge 0,\; \quad \quad   \quad  i=1,\cdots, m\;\; & \text{ in  } \Omega,\\
u^{\varepsilon}_{i}(x) =\phi_{i}(x), \;  i=1,\cdots, m\;\;    &   \text{ on   } \partial \Omega.
 \end{cases}
\end{align}
Here $\phi_{i}$ are non-negative  $C^{1}$  functions with disjoint supports that is, $\phi_{i} \cdot \phi_{j}  =0,$
almost everywhere on  the boundary, and the term $\frac{ 1 }{\varepsilon}$ is the  competition rate.

 This  model is  also called  adjacent segregation,  modeling  when  particles annihilate each other  on contact. The system (\ref{s1}) has been generalized
for  nonlinear diffusion or long segregation, where species interact at a distance from each other see \cite{4}. Also  in \cite{10} the generalization of this  problem has been considered  for the extremal Pucci operator. The numerical treatment of the limiting  case in  system (\ref{s1}) is given in \cite{2}.

The limiting configuration (solution) of (\ref{s1})  as   $\varepsilon $  tends to zero,  is related to  a free boundary problem and the  densities $u_i$  satisfy the   system of  differential inequalities. The uniqueness of  limiting solution is proven for  the cases $m=2$ in \cite{5} and $m=3$ in planar domain, see \cite{6}.  Later  in \cite{11} these  uniqueness results have  been  generalized  to arbitrary dimension and arbitrary number of species.

In this   work we give original  proof   for  uniqueness of the limiting configuration for arbitrary $m$ competing densities  by employing properties of limiting solution, which is  different approach and straightforward.
%

 The outline of the paper is as follows:  We  state the problem in Section $2$ and  provide mathematical background   and known results,  which will be used in our proof. In  Section $3$  we prove the uniqueness  of the system (\ref{s1})   in the limiting case as  $\varepsilon$ tends to zero.

\section{Known Results and Mathematical Background}

In this section we mention some of known results  for the solutions  of the system  (\ref{s1}), which  play an important role in our study. Namely, we recall some estimates  and compactness properties. 

To start with, for each fixed  $\varepsilon,$ the system (\ref{s1}) has a unique solution, see  \cite{11}.
The authors in \cite{11} use  the sub- and sup-solution method for  nonlinear elliptic systems  to construct iterative  monotone sequences  which leads to the uniqueness in case of  system (\ref{s1}).

 Let $ U^{\varepsilon}=(u_{1}^{\varepsilon},  \cdots ,u_{m}^{\varepsilon} ) $ be the unique  solution of the system (\ref{s1}) for fixed $\varepsilon.$ Then   $u_{i}^{\varepsilon} $  for $i=1, \cdots,m, $  satisfies  the following differential inequality:
 \begin{equation}\label{ave1}
 -\Delta u_{i}^{\varepsilon}  \le 0 \quad \text{ in } \quad \Omega.
  \end{equation}
Define $ \widehat{u}_{i}^{\varepsilon}$  as
 \[
  \widehat{u}_{i}^{\varepsilon}  := u_{i}^{\varepsilon}  - \sum_{j\neq i}u_{j}^{\varepsilon},
  \]
   then it is easy to verify the following property
   \begin{equation}\label{ave2}
 -  \Delta    \widehat{u}_{i}^{\varepsilon} =  \sum_{j\neq i}\sum_{h\neq j}u_{j}^{\varepsilon} \, u_{h}^{\varepsilon} \ge 0.
  \end{equation}

By constructing of sub and super solution to the system (\ref{s1}), we can show that
   $\frac{ \partial u_{i}^{\varepsilon}}{  \partial  n} $ is bounded on $ \partial \Omega$ (independent of $ \varepsilon$).
   Then multiplying the inequality $ -\Delta u_{i}^{\varepsilon}  \le 0 $   by  $ u_{i}^{\varepsilon} $    and  integrating by part yields that  $ u_{i}^{\varepsilon} $ is bounded  in  $H^{1}(\Omega)$  for each $ \varepsilon$.

The above discussion  show that the solution of  the system (\ref{s1})  belongs to the following class $F$,  see Lemma $2.1$ in  \cite{5}.
\begin{equation*}
\textit{F}=\left\{(u_1, \cdots ,u_m)  \in (H^{1}(\Omega))^{m}: u_{i} \ge 0, \,- \Delta u_{i}\leq 0, \, -\Delta  \widehat{u}_{i}\geq 0,\, u_{i}=\phi_{i}  \text{ on }    \partial \Omega    \right\},
\end{equation*}
where as in system(\ref{s1}) the boundary data  $\phi_{i}  \in C^{1}(\partial \Omega),$ nonnegative   functions and  $\phi_{i} \cdot \phi_{j}  =0,$ almost everywhere on  the boundary.

 The following result  in \cite{3,5}  shows the asymptotic behavior of the system as   $ \varepsilon    \rightarrow  0.$
     Let $U^{\varepsilon}=(u_{1}^{\varepsilon},\cdots ,u_{m}^{\varepsilon})
$ be the  solution of  system (\ref{s1})   for a fixed $\varepsilon$.  If $ \varepsilon $ tends to zero, then there exists $ U=(u_1, \cdots,u_m) \in (H^{1}(\Omega))^{m} $ such
that for all $ i=1,\cdots,m$:

\begin{enumerate}
        \item up to a  subsequences, $ u_{i}^{\varepsilon}\rightarrow u_{i} $
\  strongly in $H^{1}(\Omega)$,
        \item $ u_{i}\cdot u_{j}=0 $ if $ i\neq j$ a.e in \ $ \Omega$,
        \item $\Delta u_{i}=0 $ in the set $ {\{u_{i}>0}\}$,
        \item  Let $x$  belongs to the common interface  of two components $u_i$ and $u_j,$  then\\
        $
   \underset{ y\rightarrow x} {\text{ lim}} \ \nabla u_{i}(y)=-
\underset{ y\rightarrow x}{ \text{ lim}} \  \nabla u_{j}(y).
  $
\end{enumerate}
%
%
%
%
From above the limiting solution, as $\varepsilon $ tends to zero, belongs to the following  class:
\[
\textit{S}=\left\{(u_1, \cdots ,u_m) \in F: u_i\cdot u_j=0 \quad  \text{for } i\neq j    \right\}.
\]
Note that the inequalities in (\ref{ave1}) and (\ref{ave2}) hold as $\varepsilon$ tends to zero. Also
$$
-  \Delta    \widehat{u}_{i} = 0 \quad \text{ on  } \quad {\{x\in \Omega:\,   u_{i}(x) > 0}\}.
  $$

	In this part we briefly review the known    results  about  uniqueness of  the limiting configuration of the system (\ref{s1}).
   In particular,  for the case $m=2,$ the  limiting solution and the rate of convergence are  given (see  Theorem $2.1$ in \cite{5}).  For the sake of clarity we recall that result below.

 \begin{theorem}   Let $W$  be  harmonic  in $\Omega$  with the boundary data $\phi_{1}-\phi_{2} $. Let $u_{1}~=~W^{+}$,  $u_{2}=-
W^{-}$, then  the pair $(u_{1},u_{2})$ is the limit configuration of
any sequences  $ (u_{1}^{\varepsilon},u_{2}^{\varepsilon})$ and:
  \[
  \parallel u_{i}^{\varepsilon}-u_{i}
\parallel_{H^{1}(\Omega)}\leq C \cdot {\varepsilon}^{1/6} \ \ as \ \varepsilon\rightarrow 0, \quad i=1,2.
  \]
\end{theorem}

For the case $m=3,$ the uniqueness of the limiting configuration, as $ \varepsilon $   tends to zero, is shown in \cite{6}  on a planar domain, with appropriate boundary conditions.  More precisely, the authors prove that the limiting configuration of the following system

\begin{equation*}
\left \{
\begin{array}{lll}
 \Delta u_{i}^{\varepsilon}= \frac{u_{i}^{\varepsilon}(x)}{\varepsilon} \sum\limits_{j\neq i}^{3}   u_{j}^{\varepsilon}(x)  & \text{ in   }  \Omega,\\
u^{\varepsilon}_{i}(x) =\phi_{i}(x) \;     &   \text{ on   } \partial \Omega,\\
  i = 1, 2, 3,
\end{array}
\right.
\end{equation*}
    minimizes the energy
\begin{equation*}
  E(u_1, u_2, u_3)=\int_{\Omega}  \sum_{i=1}^{3}  \frac{1}{2}| \nabla u_{i}|^{2} dx,
\end{equation*}
among all segregated states   $u_i \cdot  u_j = 0 ,$  a.e.  with the same boundary conditions.

\begin{remark}
The system (\ref{s1}) is not in a  variational form. In \cite{7} for a class of segregation states governed by a variational
principle the proof of existence and uniqueness are shown under some suitable conditions.
\end{remark}

In  \cite{11}   the  uniqueness of the limiting configuration and least energy property  are generalized  to arbitrary dimension  and  for  arbitrary number of components. Following notations in \cite{11}, let  $ \sum $   denote the metric space

\[
 {\{ (u_1, u_2 , \cdots, u_m) \in \mathbb{R}^m :  u_{i}\ge 0,\, u_{i}u_{j}=0 \quad \text{for} \quad i\neq  j\}}.
    \]
The authors in  \cite{11}     show that  the solution of the limiting problem   $(u_1, \cdots ,u_m) \in S$ is a harmonic map into the space  $ \sum. $      The harmonic map is the critical point (in weak sense) of the following energy  functional
	\[
	\int_{\Omega}  \sum_{i=1}^{m}  \frac{1}{2}| \nabla u_{i}|^{2} dx,
	\]
	among all nonnegative  segregated states   $u_i \cdot  u_j = 0 ,$  a.e.  with the same boundary conditions,  see Theorem 1.6 in \cite{11}.

Their proof   is based on  computing  the derivative of the energy functional with respect to the geodesic homotopy between  $u  $ and a comparison  to an energy minimizing map $v$  with same boundary values. This demands some procedures  to  avoid  singularity of free boundary.  Unlike their approach,  our proof   is more direct and based on properties of limiting solutions  and doesn't require results from regularity theory or harmonic maps.

%
\section{Uniqueness }
  In this section  we prove the uniqueness for the limiting case as $\varepsilon $ tends to zero. Our approach is motivated from the recent work related to the numerical analysis of a certain class of the spatial segregation of reaction-diffusion systems (see \cite{1}). We  heavily use the following notation:
  $$
  \widehat{w}_{i}(x):={w}_{i}(x)-\sum_{p\neq i} {w}_{p}(x),
  $$
 for every $1\le i \le m.$

\begin{lemma}\label{s11}
  Let   two  elements  $(u_1, \cdots , u_{m} ) $ and
 $(v_{1}, \cdots , v_{m} ) $ belong to $S$. Then the following equation for each $ 1 \le i  \le m$ holds:
 \[
\underset{\overline{\Omega}}  {\max} (\widehat{u}_{i}(x)  - \widehat{v}_{i}(x))= \underset{{\{  u_{i}(x) \le  v_{i}(x)}\}}  {\max} (\widehat{u}_{i}(x)  - \widehat{v}_{i}(x)).
\]
\end{lemma}
\begin{proof}
We argue by contradiction. Let  there exists some  $i_0$ such that
\begin{equation}\label{contr-cond}
 \underset{\overline{\Omega}}  {\max} (\widehat{u}_{i_0} - \widehat{v}_{i_0})= \underset{{\{  u_{i_0} >   v_{i_0}}\}}  {\max} (\widehat{u}_{i_0} -\widehat{v}_{i_0})>  \underset{{\{  u_{i_0} \leq   v_{i_0}}\}}  {\max} (\widehat{u}_{i_0} -\widehat{v}_{i_0}).
\end{equation}
Assume $D={\{ x \in \Omega :  u_{i_0}(x) >   v_{i_0}(x)}\}, $  then in $D$ we have

\begin{equation}\label{s5}
\left \{
\begin{array}{lll}
-\Delta  \widehat{u}_{i_0}(x) = 0, \\
-\Delta  \widehat{v}_{i_0}(x) \ge 0,
 \end{array}
\right.
\end{equation}
which implies that
\[
\Delta  (  \widehat{u}_{i_0}(x)- \widehat{v}_{i_0}(x))   \ge  0.
\]
The weak maximum principle yields
$$
 \underset{D}  {\max} (\widehat{u}_{i_0} - \widehat{v}_{i_0})\le \underset{\partial D}  {\max} (\widehat{u}_{i_0} - \widehat{v}_{i_0}) \le \underset{{\{  u_{i_0} = v_{i_0}}\}}  {\max} (\widehat{u}_{i_0} -\widehat{v}_{i_0}),
$$
which is inconsistent with our assumption \eqref{contr-cond}. It is clear that we can interchange the role of $\widehat{u}_{i}$ and $\widehat{v}_{i}.$ Thus, we also have
\[
\underset{\overline{\Omega}}  {\max} (\widehat{v}_{i}(x)  - \widehat{u}_{i}(x))= \underset{{\{ v_{i}(x) \le  u_{i}(x)}\}}  {\max} (\widehat{v}_{i}(x)  - \widehat{u}_{i}(x)),
\]
for all $1\le i\le m.$

\end{proof}

In view of Lemma \ref{s11} we define the following quantities
 \[
 P:= \underset{1\le i\le m }{\max} \left( \underset{\overline{\Omega}}  {\max} (\widehat{u}_{i}(x)  - \widehat{v}_{i}(x))\right)=\underset{1\le i\le m }{\max} \left( \underset{ \{u_{i}\le  v_{i}\}}  {\max} (\widehat{u}_{i}(x)  - \widehat{v}_{i}(x))\right),
 \]
  \[
 Q:= \underset{1\le i\le m }{\max} \left( \underset{\overline{\Omega}}  {\max} (\widehat{v}_{i}(x)- \widehat{u}_{i}(x) ) \right)=\underset{1\le i\le m }{\max} \left( \underset{ \{v_{i}\le  u_{i}\}}  {\max} (\widehat{v}_{i}(x)  - \widehat{u}_{i}(x))\right).
 \]
\begin{lemma}\label{s21}
Let   two  elements  $(u_1, \cdots , u_{m} ) $ and $(v_{1}, \cdots , v_{m} ) $ belong to $S$. We set $P$ and $Q$ as defined above.
 If $P> 0$ is attained for some index $1\le i_0\le m,$  then  we have $P=Q> 0.$ Moreover, there exist another index $ j_0\neq i_0 $ and a point $x_0\in \Omega, $  such that:
   \[
   P=Q=\underset{{\{u_{i_0} \le  v_{i_0}}\} }  {\max} (\widehat{u}_{i_0}  - \widehat{v}_{i_0})=    \underset{{\{u_{i_0} = v_{i_0}=0}\} }  {\max} (\widehat{u}_{i_0}  - \widehat{v}_{i_0})=   v_{j_0}(x_0)  -  u_{j_0}(x_0).
   \]
\end{lemma}
\begin{proof}
Let  the maximum $P>0$ be  attained  for the ${i_0}^{\text{th}}$ component.  According to the previous lemma,  we know  that
$(\widehat{u}_{i_0}(x)  - \widehat{v}_{i_0}(x)) $  attains its maximum on the set $ {\{  u_{i_0}(x)\leq  v_{i_0}(x)}\}.$ Let that maximum point be $x^*  \in {\{  u_{i_0}(x)\leq  v_{i_0}(x)}\} .$  It is easy to see that    $\widehat{u}_{i_0}(x^*)  - \widehat{v}_{i_0}(x^*)=P>0,$ implies $ u_{i_0}(x^*)= v_{i_0}(x^*)=0.$ Indeed,  if $ u_{i_0}(x^*)= v_{i_0}(x^*)>0,$ then
 in the light of disjointness property of the components  of ${u}_{i}$ and ${v}_{i}$ we get $P=\widehat{u}_{i_0}(x^*)  - \widehat{v}_{i_0}(x^*)={u}_{i_0}(x^*)  - {v}_{i_0}(x^*)=0$ which is a contradiction. If $u_{i_0}(x^*)<v_{i_0}(x^*),$ then again due to the disjointness of the densities $u_i,v_i,$ we have
 $$
0<P=\widehat{u}_{i_0}(x^*)  - \widehat{v}_{i_0}(x^*)=\widehat{u}_{i_0}(x^*) -{v}_{i_0}(x^*) \leq {u}_{i_0}(x^*)  - {v}_{i_0}(x^*)<0.
 $$
 This again leads to a contradiction. Therefore ${u}_{i_0}(x^*)={v}_{i_0}(x^*)=0.$


Now assume by contradiction that $ Q  \le 0.$ Then  by definition of $Q$ we  should have
\[
  \widehat{v}_{j}(x)  \le  \widehat{u}_{j}(x), \quad \forall x \in \Omega ,\;  j=1, \cdots, m.
\]
This apparently yields
 \[
   v_{j}(x)  \le  u_{j}(x), \quad \forall x \in \Omega ,\,  j=1, \cdots, m.
   \]
   Let  $D_{i_0}= {\{  u_{i_0}(x)=   v_{i_0}(x)=0}\},$ then we have
   \[
   0< P= \underset { D_{i_0}} {\max}\left (\widehat{u}_{i_0}(x)  - \widehat{v}_{i_0}(x)\right)=\underset { D_{i_0}} {\max}\; \left(\sum_{j\neq {i_0}} (v_{j}(x)- u_{j}(x))\right) \le 0.
   \]
   This contradiction implies that $Q> 0$. By analogous proof, one can see that if $P$ be non-positive then $Q$  will be non-positive  as well.
   Next, assume the maximum  $P$ is attained at a point $x_0 \in D_{i_0}.$
    Then, we get
   \begin{multline}
   	 0< P=  \widehat{u}_{i_0}(x_0)  - \widehat{v}_{i_0}(x_0)= (u_{i_0}(x_0)  - v_{i_0}(x_0)) +\\+  \sum_{j\neq {i_0}}( v_{j} (x_0) - u_{j}(x_0))= \sum_{j\neq {i_0}}( v_{j} (x_0) - u_{j}(x_0)).  	
   \end{multline}
This  shows that
\[
 \sum_{j\neq i_0}   v_{j} (x_0)= \sum_{j\neq i_0}  u_{j}(x_0)+P> 0.
\]
Since $(v_1, \cdots ,v_m)\in S, $ then there exists $j_0\neq i_0$ such that $v_{j_0}(x_0)> 0. $ This implies
\[
 0< P=  \widehat{u}_{i_0}(x_0)  - \widehat{v}_{i_0}(x_0)= v_{j_0}(x_0)- \sum_{j\neq i_0}u_{j}(x_0) \le \widehat{v}_{j_0}(x_0)-\widehat{u}_{j_0}(x_0) \le Q.
 \]
 The same argument shows  that $Q \le P$ which yields $P=Q$. Hence, we can write
 \[
 P=v_{j_0}(x_0)-\sum_{j\neq i_0}  u_{j}(x_0)=    \widehat{v}_{j_0}(x_0)-\widehat{u}_{j_0}(x_0)=Q.
  \]
This gives us  $2 \sum_{j\neq j_0}u_{j}(x_0)= 0,$ and  therefore
\[
u_{j} (x_0)=0,   \quad   \forall j\neq j_0,
\]
which completes the  last statement of the  proof.
\end{proof}
We are ready to prove the uniqueness of a limiting configuration.
\begin{theorem}
  There exists a unique vector  $(u_1,\cdots,u_m) \in S ,$  which satisfies the limiting solution of (\ref{s1}).
\end{theorem}

\begin{proof}
In order to show the uniqueness of the limiting configuration,  we assume that  two m-tuples    $(u_1,\cdots,u_{m})$ and $(v_{1},\cdots,v_{m} )$ are  the solutions of  the system  (\ref{s1})  as
$ \varepsilon$ tends to zero.   These  two solutions belong to the class  $S$. For them  we set $P$ and $Q$ as above.
  Then, we consider two cases  $P\le  0 $  and $P> 0.$
If we assume that $ P \le  0$  then     Lemma \ref{s21} implies that  $Q \le  0$.  This leads to
\[
0 \le -Q\le \widehat{u}_{i}(x)  - \widehat{v}_{i}(x) \le P\le 0,
\]
for every $1\le i\le m,$ and $x\in\Omega.$
This provides that
\[
 \widehat{u}_{i}(x)  =\widehat{v}_{i}(x) \quad  i= 1,\cdots ,m,
\]
  which in turn implies
\[
u_{i}(x) =v_{i}(x).
\]

Now, suppose  $P > 0.$  We show that this  case leads to a contradiction.
Let the value $P$ is attained for some  $i_0,$  then due to Lemma \ref{s21} there exist $ x_0 \in \Omega  $  and
  $ j_0\neq  i_0 $ such that:
  \[
  0 < P= Q = \widehat{u}_{i_0}(x_0) - \widehat{v}_{i_0}(y_0)=   \underset{{\{ u_{i_0}=v_{i_0}=0 }\}} {\max}(\widehat{u}_{i_0}(x) - \widehat{v}_{i_0}(x))= v_{j_0}(x_0)- u_{j_0}(x_0).
  \]
 	Let $\Gamma$ be a fixed curve starting at  $x_0$ and ending on  the boundary of $\Omega.$  Since $\Omega$ is connected, then one can always choose such a curve belonging to $\overline{\Omega}.$  	By the disjointness and smoothness  of  $v_{j_0}$ and  $u_{j_0}$ there exists a  ball centered at $x_0$, and with radius  $r_{0}$ ($ r_{0} $  depends on $x_{0}$)  which  we denote it $B_{r_{0}}(x_0)$,    such that
 \[
 {v}_{j_0}(x)-u_{j_0}(x)>0 \text{ in } \quad B_{r_{0}}(x_0).
 \]
 This yields
 	\[
 	\Delta (\widehat{v}_{j_0}(x)-\widehat{u}_{j_0}(x))\ge 0\;\;\mbox{in}\;\;  B_{r_{0}}(x_0).
 	\]
 	The maximum principle implies that
 	$$
 \underset{ \overline{B_{r_{0}}(x_0)}}{\max}\;(\widehat{v}_{j_0}(x)-\widehat{u}_{j_0}(x))= \max_{\partial  B_{r_{0}}(x_0)}\;(\widehat{v}_{j_0}(x)-\widehat{u}_{j_0}(x))\leq P.
 	$$
 	One the other hand, in view of  Lemma \ref{s21} we have
 	$$
    \widehat{v}_{j_0}(x_0)-\widehat{u}_{j_0}(x_0)=v_{j_0}(x_0)- u_{j_0}(x_0)=P,
 	$$
 	which implies that $P$ is attained at the interior point $x_0\in  B_{r_0}(x_0).$ Thus,
 	$$
 	\widehat{v}_{j_0}(x)-\widehat{u}_{j_0}(x)\equiv P>0\;\;\mbox{in}\;\; \overline{B_{r_{0}}(x_0)}.
 	$$
 	Next  let $x_1\in  \Gamma\cap \partial B_{r_{0}}(x_0).$  We get
 	$\widehat{v}_{j_0}(x_1)-\widehat{u}_{j_0}(x_1)=P>0,$ which leads to ${v}_{j_0}(x_1)~\ge ~ {u}_{j_0}(x_1).$ We proceed as follows: If ${v}_{j_0}(x_1)> {u}_{j_0}(x_1),$ then as above ${v}_{j_0}(x)> {u}_{j_0}(x)$ in $B_{r_{1}}(x_1).$ This in turn implies
 	\[
 	\Delta (\widehat{v}_{j_0}(x)-\widehat{u}_{j_0}(x))\ge 0\;\;\mbox{in}\;\;  B_{r_{1}}(x_1).
 	\]
 	Again following  the maximum principle and recalling that $\widehat{v}_{j_0}(x_1)-\widehat{u}_{j_0}(x_1)=P$ we conclude that
 	$$
 	\widehat{v}_{j_0}(x)-\widehat{u}_{j_0}(x)= P>0\;\;\mbox{in}\;\; \overline{B_{r_{1}}(x_1)}.
 	$$
 		
 	If ${v}_{j_0}(x_1)={u}_{j_0}(x_1),$ then clearly the only possibility is ${v}_{j_0}(x_1)={u}_{j_0}(x_1)=0.$ Thus,
 	$$
 	0<P=\widehat{v}_{j_0}(x_1)-\widehat{u}_{j_0}(x_1)=\sum_{j\neq j_0}({u}_{j}(x_1)-{v}_{j}(x_1)).
 	$$
 	Following the lines of the proof of Lemma \ref{s21}, we find some $k_0\neq j_0,$ such that
 	$$
 	P={u}_{k_0}(x_1)-{v}_{k_0}(x_1)=\widehat{u}_{k_0}(x_1)-\widehat{v}_{k_0}(x_1).
 	$$
 	 	It is easy to see that there exists a ball  $B_{r_{1}}(x_1)$ (without loss of generality one keeps the same notation)
 	\[
 	\Delta (\widehat{u}_{k_0}(x)-\widehat{v}_{k_0}(x))\ge 0\;\;\mbox{in}\;\;  B_{r_{1}}(x_1).
 	\]
 	In view of the maximum principle and above steps we obtain:
 	$$
 	\widehat{u}_{k_0}(x)-\widehat{v}_{k_0}(x)= P>0\;\;\mbox{in}\;\; \overline{B_{r_{1}}(x_1)}.
 	$$
 	
 	Then we take $x_2\in \Gamma\cap \partial B_{r_{1}}(x_1)$  such that $x_1$ stands between the points $x_0$ and $x_2$ along the given curve $\Gamma.$ According to the previous arguments for the point $x_2$ we will find an index $l_0\in{\{1,\cdots,m}\}$ and corresponding ball $ B_{r_{2}}(x_2),$ such that
 	$$
 	|\widehat{u}_{l_0}(x)-\widehat{v}_{l_0}(x)|=P \;\;\mbox{in}\;\; \overline{B_{r_{2}}(x_2)}.
 	$$
 	We continue this way and obtain a sequence of  points $x_n$   along the curve $\Gamma$, which are getting closer to the boundary of $\Omega.$ Since for all $j=1, \cdots ,m$ and $x\in\partial\Omega$ we have
 	$$
 	\widehat{u}_{j}(x)-\widehat{v}_{j}(x)=\widehat{v}_{j}(x)-\widehat{u}_{j}(x)=0,
 	$$
 	then obviously after finite steps $N$ we find the point $x_N,$ which will be very close  to the $\partial\Omega$ and for all $j=1,\cdots,m $
 	$$
 	|\widehat{u}_{j}(x_N)-\widehat{v}_{j}(x_N)|<P/2.
 	$$
 	On the other hand, according to our construction for the point $x_N,$ there exists an index $1\le j_N \le m$ such that
 	$$
 	|\widehat{u}_{j_N}(x_N)-\widehat{v}_{j_N}(x_N)|=P,
 	$$
 	which is a contradiction. This completes the proof of the uniqueness.

\end{proof}

\bibliographystyle{acm}%
\bibliography{seg-unique}

\end{document}